\documentclass[reqno,11pt]{article}

\usepackage{amsmath,latexsym,amssymb}

\def\sqr#1#2{{\vcenter{\hrule height.#2pt
        \hbox{\vrule width.#2pt height#1pt \kern#1pt
                \vrule width.#2pt}
        \hrule height.#2pt}}}
\def\square{\mathchoice\sqr64\sqr64\sqr{4}3\sqr{3}3}
\def\QED{\hfill$\square$\break}
\def\demo{\noindent{\bf Proof: }}

\newtheorem{Theorem}{\sc Theorem}[section]
\newtheorem{Lemma}[Theorem]{\sc Lemma}
\newtheorem{Corollary}[Theorem]{\sc Corollary}
\newtheorem{Proposition}[Theorem]{\sc Proposition}

\newtheorem{Remark}[Theorem]{\sc Remark}

\newtheorem{Definition}[Theorem]{\sc Definition}

\DeclareMathOperator{\content}{c}

\DeclareMathOperator{\V}{V}
\def\Z{\mathbb Z}
\def\N{\mathbb N}
\DeclareMathOperator{\rank}{rank}
\DeclareMathOperator{\Coker}{Coker} 
\DeclareMathOperator{\HH}{H}
\DeclareMathOperator{\depth}{depth}

\begin{document}

\baselineskip=13pt

\pagestyle{empty}

\ \vspace{1.7in}

\noindent {\LARGE\bf The support of top graded local cohomology modules}

\vspace{.25in}

\noindent Mordechai Katzman, \  Department of Pure Mathematics,
University of Sheffield, United Kingdom.
{\it E-mail}: {\tt M.Katzman@sheffield.ac.uk}

\vspace{2.4cm}

\section{Introduction \hfill\break}

Let $R_0$ be any domain,
let $R=R_0[U_1, \dots, U_s]/I$, where $U_1, \dots, U_s$
are indeterminates of positive degrees $d_1,\dots, d_s$, and
$I\subset R_0[U_1, \dots, U_s]$ is a homogeneous ideal.

The main theorem in this paper is Theorem \ref{theorem1},
a generalization of Theorem 1.5 in \cite{KS},
which states that
all the associated primes of
$H:=H^s_{R_+}(R)$ contain a certain non-zero ideal $\content(I)$ of $R_0$ called the
``content'' of $I$ (see Definition \ref{definition1}.)
It follows that the support of $H$ is simply $\V(\content(I)R + R_+)$ (Corollary 1.8) and, in particular,
$H$ vanishes if and only if $\content(I)$ is the unit ideal.

These results raise the question of whether local cohomology modules have finitely many minimal associated primes--
this paper provides further evidence in favour of such a result (Theorem \ref{theorem2} and Remark \ref{remark1}.)

Finally, we give a very short proof of a weak version of the monomial conjecture based on Theorem \ref{theorem1}.

\bigskip
\section{\bf The vanishing of top local cohomology modules\hfill\break}\label{section1}

Throughout this section $R_0$ will denote an arbitrary commutative Noetherian domain.
We set $S=R_0[U_1, \dots, U_s]$ where $U_1,\dots, U_s$ are indeterminates of degrees $d_1, \dots, d_s$,
and $R=S/I$ where $I\subset R_0[U_1, \dots, U_s]$ is an homogeneous ideal.
We define $\Delta=d_1 + \dots + d_s$ and
denote with $\mathcal{D}$ the sub-semi-group of $\mathbb{N}$ generated by $d_1, \dots, d_s$.

For $t \in \Z$, we shall denote by $(\: {\scriptscriptstyle \bullet} \:)(t)$
the $t$-th shift functor (on the category of graded $R$-modules and homogeneous homomorphisms).

For any multi-index $\lambda=(\lambda^{(1)}, \dots, \lambda^{(s)})
\in \Z^s$ we shall write $U^\lambda$ for
$U_1^{\lambda^{(1)}} \dots U_s^{\lambda^{(s)}}$ and we shall set
$| \lambda | = \lambda^{(1)}+ \dots+ {\lambda^{(s)}}$.

\begin{Lemma}\label{lemma1}
Let $I$ be generated by homogeneous elements
$f_1, \dots, f_r\in S$.
Then there is an exact sequence of graded $S$-modules and
homogeneous homomorphisms
$$ \bigoplus_{i=1}^r H^s_{S_+}(S)(-\deg f_i)\xrightarrow[]{(f_1, \dots, f_r)}
H^s_{S_+}(S) \longrightarrow H^s_{R_+}(R) \longrightarrow 0.$$
\end{Lemma}

\demo
The functor $H^s_{S_+}(\bullet)$ is right exact and the natural equivalence
between $H^s_{S_+}(\bullet)$ and
$(\: {\bullet} \:)\otimes_S H^s_{S_+}(S)$ (see \cite[6.1.8 \& 6.1.9]{BS})
actually yields a homogeneous $S$-isomorphism
$$H^s_{S_+}(S)/(f_1,\ldots,f_r)H^s_{S_+}(S) \cong H^s_{S_+}(R).$$ To complete the
proof, just note that there is an isomorphism of graded
$S$-modules $H^s_{S_+}(R) \cong H^s_{R_+}(R)$, by the Graded
Independence Theorem \cite[13.1.6]{BS}.
\QED

We can realize $H^s_{S_+}(S)$ as the module
$R_0[U_1^-,\dots,U_s^-]$ of inverse polynomials described in \cite[12.4.1]{BS}:
this graded $R$-module vanishes beyond degree $-\Delta$.
More generally
$R_0[U_1^-,\dots,U_s^-]_{-d} \neq 0$ if and only if $d\in \mathcal{D}$.

For each
$d \in \mathcal{D}$, $R_0[U_1^-,\dots,U_s^-]_{-d}$
is a free $R_0$-module with base
${\mathcal{B}(d)}:=\left(U^\lambda\right)_{-\lambda\in \N^s, |\lambda| = -d}. $
We combine this realisation with the previous
lemma to find a presentation of each homogeneous component of
$H^s_{R_+}(R)$ as the cokernel of a matrix with entries in $R_0$.

Assume first that $I$ is generated by one homogeneous element $f$
of degree $\delta$. For any $d\in\mathcal{D}$ we have, in view of Lemma
\ref{lemma1}, a graded exact sequence
$$R_0[U_1^-,\dots, U_s^-]_{-d-\delta} \xrightarrow[]{\phi_d}
R_0[U_1^-,\dots, U_s^-]_{-d} \longrightarrow H^s_{R_+}(R)_{-d}
\longrightarrow 0 .$$ The map of free $R_0$-modules $\phi_d$ is given
by multiplication on the left by a $\#{\mathcal{B}(d)} \times \#{\mathcal{B}(d+\delta)}$ matrix which we shall denote later
by $M(f;d)$.

In the general case, where $I$ is generated by homogeneous
elements $f_1, \dots, f_r\in S$, it follows from Lemma
\ref{lemma1} that the $R_0$-module $H^s_{R_+}(R)_{-d}$ is the
cokernel of a matrix $M(f_1, \dots, f_r; d)$ whose columns consist
of all the columns of $M(f_1,d), \dots, M(f_r,d)$.

Consider a homogeneous $f\in S$ of degree $\delta$. We shall now
describe the matrix $M(f;d)$ in more detail and to do so we start
by ordering the bases of the source and target of $\phi_d$ as
follows. For any $\lambda, \mu\in \Z^s$ with negative entries we
declare that $U^\lambda < U^\mu$ if and only if $U^{-\lambda}
<_{\mathrm{Lex}} U^{-\mu}$ where ``$<_{\mathrm{Lex}}$'' is the
lexicographical term ordering in $S$ with $U_1 > \dots > U_s$.
We order the bases ${\mathcal{B}(d)}$, and by doing so also the columns and rows of
$M(f;d)$, in ascending order.
We notice that the entry in $M(f;d)$ in the $U^\alpha$ row and $U^\beta$ column is now
the coefficient of $U^\alpha$ in $f U^\beta$.

\begin{Lemma}\label{lemma1.5}
Let $\nu \in \Z^s$ have negative entries and let $\lambda_1, \lambda_2\in \mathbb{N}^s$.
If $U^{\lambda_1} <_{\mathrm{Lex}}  U^{\lambda_2}$ and
$U^\nu U^{\lambda_1}$, $U^\nu U^{\lambda_2}\in R_0[U_1^-,\dots,U_s^-]$ do not vanish then
$U^\nu U^{\lambda_1} >  U^\nu U^{\lambda_2}$.
\end{Lemma}

\demo
Let $j$ be the first coordinate in which $\lambda_1$ and $\lambda_2$ differ.
Then
$\lambda_1^{(j)}<\lambda_2^{(j)}$ and so also
$-\nu^{(j)}-\lambda_1^{(j)}>-\nu^{(j)}-\lambda_2^{(j)}$; this implies that
$U^{-\nu-\lambda_1} >_{\mathrm{Lex}} U^{-\nu-\lambda_2}$ and
$U^{\nu+\lambda_1} > U^{\nu+\lambda_2}$.
\QED

\begin{Lemma}\label{lemma2}
Let $f\neq 0$ be a homogeneous element in $S$. Then, for all
$d\in \mathcal{D}$, the matrix $M(f;d)$ has maximal rank.
\end{Lemma}

\demo
We prove the lemma by producing a non-zero maximal minor of
$M(f;d)$. Write $f=\sum_{\lambda\in \Lambda} a_\lambda U^\lambda$
where $a_\lambda\in R_0\setminus\{ 0 \}$ for all $\lambda\in
\Lambda$ and let $\lambda_0$ be such that $U^{\lambda_0}$ is the
minimal member of $\left\{U^{\lambda} : \lambda \in
\Lambda\right\}$ with respect to the lexicographical term order in $S$.

Let $\delta$ be the degree of $f$. Each column of $M(f;d)$
corresponds to a monomial $U^\lambda\in {\mathcal{B}(d+\delta)}$;
its $\rho$-th entry is the coefficient of  $U^\rho$ in $f
U^\lambda \in R_0[U_1^-,\dots, U_s^-]_{-d}$.

Fix any $U^{\nu}\in {\mathcal{B}(d)}$ and consider the column
$c_\nu$ corresponding to
$U^{\nu-\lambda_0 }\in {\mathcal{B}(d+\delta)}$. The $\nu$-th entry of $c_\nu$ is
obviously $a_{\lambda_0}$.

By the previous lemma all entries in $c_\nu$ below the $\nu$th row vanish.
Consider the square submatrix of $M(f;d)$
whose columns are the $c_\nu~(\nu\in {\mathcal{B}(d)})$; its
determinant is clearly a power of $a_{\lambda_0}$ and hence is
non-zero.
\QED

\begin{Definition}\label{definition1}
For any $f \in R_0[U_1, \dots, U_s]$ write $f=\sum_{\lambda\in\Lambda} a_\lambda U^\lambda$ where $a_\lambda\in R_0$ for all
$\lambda\in \Lambda$. For such an $f \in R_0[U_1, \dots, U_s]$ we
define the {\em content\/} $\content(f)$ of $f$ to be the ideal
$\langle a_\lambda : \lambda\in \Lambda \rangle$ of $R_0$
generated by all the coefficients of $f$. If $J\subset R_0[U_1,\dots, U_s]$ is an ideal, we define its {\em content\/}
$\content(J)$ to be the ideal of $R_0$ generated by the contents
of all the elements of $J$. It is easy to see that if $J$ is
generated by $f_1, \dots, f_r$, then $\content(J) =\content(f_1)+\dots +\content(f_r)$.
\end{Definition}

\begin{Lemma}\label{lemma3}
Suppose that $I$ is generated by homogeneous elements\goodbreak
$f_1, \dots, f_r\in S$. Fix any $d\in \mathcal{D}$. Let $t:=\rank M(f_1, \dots, f_r; d)$ and let $I_d$ be the ideal generated by all
$t\times t$ minors of $M(f_1, \dots, f_r; d)$.
Then $\content(I) \subseteq \sqrt{I_d}$.
\end{Lemma}

\demo
It is enough to prove the lemma when $r=1$; let $f=f_1$. Write
$f=\sum_{\lambda\in \Lambda} a_\lambda U^\lambda$ where
$a_\lambda\in R_0\setminus\{ 0 \}$ for all $\lambda\in \Lambda$.
Assume that $\content(I) \not\subseteq \sqrt{I_d}$ and pick
$\lambda_0$ so that $U^{\lambda_0}$ is the minimal element in
$\left\{U^{\lambda} : \lambda \in \Lambda\right\}$ (with respect
to the lexicographical term order in $S$) for which
$a_\lambda \notin\sqrt{I_d}$.
Notice that the proof of Lemma \ref{lemma2} shows
that $U^{\lambda_0}$ cannot be the minimal element of
$\left\{U^{\lambda} : \lambda \in \Lambda\right\}$.

Fix any $U^{\nu}\in {\mathcal{B}(d)}$ and consider the column
$c_\nu$ corresponding to $U^{\nu-\lambda_0} \in {\mathcal{B}(d+\delta)}$.
The $\nu$-th entry of $c_\nu$ is
obviously $a_{\lambda_0}$. Lemma \ref{lemma1.5} shows that, for any other
$\lambda_1\in \Lambda$ with $U^{\lambda_1} >_{\mathrm{Lex}} U^{\lambda_0}$,
either $\nu-\lambda_0+\lambda_1$ has a non-negative entry, in which
case the corresponding term of $f U^{\nu-\lambda_0}\in R_0[U_1^-,\dots, U_s^-]_{-d}$ is zero, or
$U^\nu > U^{\nu-\lambda_0+\lambda_1}$.

Similarly, for any other $\lambda_1 \in \Lambda$ with
$U^{\lambda_1} <_{\mathrm{Lex}} U^{\lambda_0}$, either
$\nu-\lambda_0+\lambda_1$ has a non-negative entry, in which case the
corresponding term of
$f U^{\nu-\lambda_0}\in R_0[U_1^-,\dots,U_s^-]_{-d}$ is zero, or $U^\nu < U^{\nu-\lambda_0+\lambda_1}$.

We have shown that all the entries
below
the $\nu$-th row of
$c_\nu$ are in $\sqrt{I_d}$. Consider the matrix $M$ whose columns
are $c_\nu~(\nu \in {\mathcal{B}(d)})$ and let
$\overline{\phantom{X} } : R_0 \rightarrow R_0 / \sqrt{I_d}$
denote the quotient map. We have
$$0=\overline{\det(M)}=\det(\overline{M})=\overline{a_{\lambda_0}}^{{d-1} \choose s-1}$$
and, therefore, $a_{\lambda_0} \in \sqrt{I_d}$, a contradiction.
\QED

\begin{Theorem}\label{theorem1}
Suppose that $I$ is generated by homogeneous elements $f_1, \dots,
f_r\in S$. Fix any $d\in \mathcal{D}$. Then each associated prime of
$H^s_{R_+}(R)_{-d}$
contains $\content(I)$. In particular
$H^s_{R_+}(R)_{-d}=0$ if and only if $\content(I)=R_0$.
\end{Theorem}

\demo
Recall that for any $p,q \in \N$ with $p\leq q$ and any $p\times
q$ matrix $M$ of maximal rank with entries in any domain, $\Coker
M = 0$ if and only if the ideal generated by the maximal minors of
$M$ is the unit ideal. Let $M=M(f_1, \dots, f_r; d)$, so that
$H^s_{R_+}(R)_{-d}\cong\Coker M$.

In view of Lemmas \ref{lemma2} and \ref{lemma3}, the ideal
$\content(I)$ is contained in the radical of the ideal generated
by the maximal minors of $M$; therefore, for each $x\in
\content(I)$, the localization of $\Coker M$ at $x$ is zero;
we deduce that
$\content(I)$ is contained in all associated primes of $\Coker M$.

To prove the second statement, assume first that $\content(I)$ is
not the unit ideal. Since all minors of $M$ are contained in
$\content(I)$, these cannot generate the unit ideal and $\Coker M
\neq 0$. If, on the other hand,  $\content(I)=R_0$ then $\Coker M$
has no associated prime and $\Coker M = 0$.
\QED

\begin{Corollary}\label{cor0} Let the situation be as in\/ {\rm \ref{theorem1}}. The
following statements are equivalent:
\begin{enumerate}
\item $\content(I)=R_0$;
\item $H^s_{R_+}(R)_{-d}=0$ for some $d \in \mathcal{D}$;
\item $H^s_{R_+}(R)_{-d}=0$ for all $d \in \mathcal{D}$.
\end{enumerate}
Consequently, $H^s_{R_+}(R)$ is asymptotically gap-free in the sense of\/
{\rm \cite[(4.1)]{BH}}.
\end{Corollary}

\begin{Corollary}\label{cor2}
The $R$-module $H^s_{R_+}(R)$ has finitely many minimal associated
primes, and these are just the minimal primes of the ideal
$\content(I)R + R_+$.
\end{Corollary}

\demo
Let $r \in \content(I)$. By Theorem \ref{theorem1}, the
localization of $H^s_{R_+}(R)$ at $r$ is zero. Hence each
associated prime of $H^s_{R_+}(R)$ contains $\content(I)R$.  Such
an associated prime must contain $R_+$, since $H^s_{R_+}(R)$ is
$R_+$-torsion.

On the other hand, $H^s_{R_+}(R)_{-\Delta} \cong R_0/\content(I)$ and
$H^s_{R_+}(R)_i = 0$ for all $i > -\Delta$; therefore there is an
element of the $(-\Delta)$-th component of $H^s_{R_+}(R)$ that has
annihilator (over $R$) equal to $\content(I)R + R_+$. All the
claims now follow from these observations.
\QED

\begin{Remark}\label{rmk1}\rm
In \cite[Conjecture 5.1]{Hu}, Craig Huneke
conjectured that every local cohomology module (with respect to any
ideal) of a finitely generated module over a local Noetherian ring
has only finitely many associated primes.
This conjecture was shown to be false (cf. \cite[Corollary 1.3]{K})
but Corollary \ref{cor2} provides some evidence in support of the weaker
conjecture that every local cohomology module (with respect to any
ideal) of a finitely generated module over a local Noetherian ring
has only finitely many {\em minimal\/} associated primes.

The following theorem due to Gennady Lyubeznik (\cite{L}) gives further support for this conjecture:
\end{Remark}

\begin{Theorem}\label{theorem2}
Let $R$ be any Noetherian ring of prime characteristic $p$ and let $I\subset R$ be any ideal generated
by $f_1, \dots, f_s \in R$.
The support of $H^s_I(R)$ is Zariski closed.
\end{Theorem}

\demo
We first notice that the localization of $H^s_I(R)$ at a prime $P\subset R$ vanishes
if and only if there exist positive integers $\alpha$ and $\beta$ such that
$$ \left(f_1\cdot \dots \cdot f_s \right)^\alpha \in \langle f_1^{\alpha+\beta},\dots,  f_s^{\alpha+\beta} \rangle $$
in the localization $R_P$.
This is because if we can find such $\alpha$ and $\beta$ we can then take $q:=p^e$ powers and obtain
$$ \left(f_1\cdot \dots \cdot f_s \right)^{q \alpha} \in \langle f_1^{q\alpha+q\beta},\dots,  f_s^{q\alpha+q\beta} \rangle $$
for all such $q$.
This shows that all elements in the direct limit sequence
$$ R/\langle f_1, \dots f_s \rangle \xrightarrow[]{f_1\cdot \ldots \cdot f_s}
R/\langle f_1^2, \dots f_s^2 \rangle  \xrightarrow[]{f_1\cdot \ldots \cdot f_s} \dots $$
map to $0$ in the direct limit and hence $H^s_I(R)=0$.

But if
$$ \left(f_1\cdot \dots \cdot f_s \right)^\alpha \in \langle f_1^{\alpha+\beta},\dots,  f_s^{\alpha+\beta} \rangle $$
in $R_P$, we may clear denominators and deduce that this occurs on a Zariski open subset containing $P$.

Thus the complement of the support is a Zariski open subset.
\QED

\bigskip
It may be reasonable to expect that non-top local cohomology modules might also have finitely many minimal associated primes;
the only examples known to me of non-top local cohomology modules with infinitely many associated primes are the following:
Let $k$ be any field, let $R_0=k[x,y,s,t]$ and let
$S$ be the localisation of $R_0[u,v,a_1,\dots,a_n]$ at the maximal ideal $\frak m$ generated by
$x,y,s,t,u,v,a_1,\dots,a_n$.  Let $f=sx^2v^2 - (t+s) xy uv + t y^2 u^2\in S$ and let $R=S/fS$.
Denote by $I$ the ideal of $S$ generated by $u,v$ and by $A$ the ideal of $S$ generated by $a_1, \dots, a_n$.

\bigskip
\begin{Theorem}
Assume that $n\geq 2$. The local cohomology module $\HH^2_{I\cap A} (R)$ has infinitely many associated primes
and $\HH^{n+1}_{I\cap A}(R)\neq 0$.
\end{Theorem}

\demo
Consider the following segment of the Mayer-Vietoris sequence
$$ \dots \rightarrow \HH^2_{I+A}(R) \rightarrow \HH^2_{I}(R) \oplus \HH^2_{A}(R) \rightarrow \HH^2_{I\cap A}(R) \rightarrow \dots $$

Notice that $a_1,\dots,a_n,u$ form a regular sequence on $R$ so $\depth_{I+A} R\geq n+1\geq 3$ and the leftmost module vanishes. Thus
$\HH^2_{I}(R)$ injects into $\HH^2_{I\cap A}(R)$ and Corollary 1.3 in \cite{K} shows that $\HH^2_{I\cap A}(R)$ has infinitely many associated primes.

Consider now the following segment of the Mayer-Vietoris sequence
$$ \dots \rightarrow \HH^{n+1}_{I\cap A}(R) \rightarrow \HH^{n+2}_{I+A}(R) \rightarrow  \HH^{n+2}_{I}(R) \oplus \HH^{n+2}_{A}(R) \rightarrow \dots $$
The direct summands in the rightmost module vanish since both $I$ and $A$ can be generated by less than $n+2$ elements, so
$\HH^{n+1}_{I\cap A}(R)$ surjects onto $\HH^{n+2}_{I+A}(R)$.

Now $\content(f)$ is the ideal of $R_0$ generated by $sx^2, -(t+s)xy$ and $ty^2$ so
$\content(f)\subset \langle x,y \rangle \neq R_0$.
Corollary \ref{cor0} now shows that $\HH^{n+2}_{I+A}(R)$ does not vanish and, therefore, nor does $\HH^{n+1}_{I\cap A}(R)$.
\QED

\begin{Remark}\label{remark1}\rm
When $n\geq 3$, $\HH^3_{I+A}(R)=0$ and the argument above shows that
$\HH^2_{I}(R) \oplus \HH^2_{A}(R) \cong \HH^2_{I\cap A}(R)$.
Corollary \ref{cor2} implies that $\HH^2_{I}(R)$ has finitely many minimal primes and since the
only associated prime of $\HH^2_{A}(R)$ is $A$, $\HH^2_{I\cap A}(R)$ has finitely many minimal primes.

When $n=2$
we obtain a short exact sequence
$$ 0 \rightarrow \HH^2_{I}(R) \oplus \HH^2_{A}(R) \rightarrow \HH^2_{I\cap A}(R)  \rightarrow \HH^3_{I+A}(R) \rightarrow 0.$$
The short exact sequence
$$ 0 \rightarrow S \stackrel{f}{\rightarrow} S \rightarrow R \rightarrow 0$$
implies that $\HH^3_{I+A}(R)$ injects into the local cohomology module $\HH^4_{I+A}(S)$ whose only associated prime is $I+A$, so again we see that
$\HH^2_{I\cap A}(R)$ has finitely many minimal associated primes.
\end{Remark}

\bigskip
\section{An application: a weak form of the Monomial Conjecture.\hfill\break}

In \cite{Ho} Mel Hochster suggested reducing the Monomial Conjecture to the problem of showing the vanishing
of certain local cohomology modules which we now describe.

Let $C$ be either $\mathbb{Z}$ or a field of characteristic $p>0$,
let $R_0=C[A_1, \dots, A_s]$ where $A_1, \dots, A_s$ are indeterminates,
$S=R_0[U_s, \dots, U_s]$ where $U_1, \dots, U_s$ are indeterminates
and $R=S/F_{s,t}S$ where
$$F_{s,t}=\left(U_1\cdot \ldots \cdot U_s\right)^t- \sum_{i=1}^s A_i U_i^{t+1} .$$

Suppose that
$$H_{s,t}:=H^s_{\langle U_1, \dots, U_s \rangle} (R)$$
vanishes with $C=\mathbb{Z}$.
If for some local ring $T$ we can find a system of parameters $x_1, \dots , x_s$ so that
$\left(x_1 \cdot \ldots \cdot x_s\right)^t \in \langle x_1^{t+1}, \dots, x_s^{t+1}\rangle$, i.e.,
if there exist $a_1, \dots, a_s\in T$ so that
$\left(x_1 \cdot \ldots \cdot x_s\right)^t =\sum_{i=1}^t a_i x_i^{t+1}$
we can define an homomorphism $R\rightarrow T$ by mapping $A_i$ to $a_i$ and $U_i$ to $x_i$.
We can view $T$ as an $R$-module and we have an isomorphism of $T$-modules
$$H^s_{\langle x_1, \dots, x_s\rangle} (T) \cong H^s_{\langle U_1, \dots, U_s \rangle} (R) \otimes_R T$$
and we deduce that
$$H^s_{\langle x_1, \dots, x_s\rangle} (T)=0$$
but this cannot happen since $x_1, \dots, x_s$ form a system of parameters in $T$.

We have just shown that the vanishing of $H_{s,t}$ for all $t\geq 1$ implies the Monomial Conjecture in dimension $s$.
In \cite{Ho} Mel Hochster proved that these local cohomology modules vanish when $s=2$ or when $C$ has characteristic $p>0$, but
in \cite{R} Paul Roberts showed that, when $C=\mathbb{Z}$, $H_{3,2}\neq 0$, showing that Hochster's approach cannot be used for proving
the Monomial Conjecture in dimension 3. This can be generalized further:

\begin{Proposition}
When $C=\mathbb{Z}$, $H_{s,2}\neq 0$ for all $s\geq 3$.
\end{Proposition}

\demo
We proceed by induction on $s$; the case $s=3$ is proved in \cite{R}.

Assume that for some $s\geq 1$, $\alpha\geq 0$ and $\delta>\alpha$
the monomial $x_1^{\alpha}  \dots x_{s+1}^{\alpha}$ is in the ideal of
$C[x_1,\dots,x_{s+1}, a_1, \dots, a_{s+1}]$ generated by $x_1^{\alpha+\beta}, \dots, x_{s+1}^{\alpha+\beta}$ and $F_{s+1,t}$.

Define $G_{s+1, 2}$ to be the result of substituting $a_{s+1}=0$ in $F_{s+1, 2}$, i.e.,
$$G_{s+1, 2}=\left( x_1 \dots x_{s+1} \right)^2 - \sum_{i=1}^s a_i x_i^{3} .$$
If
\begin{equation}\label{eq1}
x_1^{\alpha}  \dots x_{s+1}^{\alpha} \in \langle x_1^{\alpha+\beta}, \dots, x_{s+1}^{\alpha+\beta}, F_{s+1,2}\rangle
\end{equation}
then by setting $a_{s+1}=0$ we see that
$$x_1^{\alpha} \dots x_{s+1}^{\alpha} \in \langle x_1^{\alpha+\beta}, \dots, x_{s+1}^{\alpha+\beta}, G_{s+1,2}\rangle .$$
If we assign degree 0 to $x_1, \dots, x_s$, degree 1 to $x_{s+1}$ and degree $2$ to $a_1, \dots, a_s$
then the ideal $\langle x_1^{\alpha+\beta}, \dots, x_{s+1}^{\alpha+\beta}, G_{s+1,2}\rangle$ is homogeneous and we must have
$$x_1^{\alpha}  \dots x_{s+1}^{\alpha}\in \langle x_1^{\alpha+\beta}, \dots, x_{s}^{\alpha+\beta}, G_{s+1,2}\rangle .$$
If we now set $x_{s+1}=1$ we obtain
\begin{equation}\label{eq2}
x_1^{\alpha} \dots x_{s}^{\alpha} \in \langle x_1^{\alpha+\beta}, \dots, x_{s}^{\alpha+\beta}, F_{s,2}\rangle .
\end{equation}

Now $H_{s+1,2}= 0$ if and only if for each $\beta\geq 1$ we can find an $\alpha\geq 0$ so that equation (\ref{eq1}) holds and this
implies that for each $\beta\geq 1$ we can find an $\alpha\geq 0$ so that equation (\ref{eq2}) holds which is equivalent to
$H_{s,2}= 0$. The induction hypothesis implies that $H_{s,2}\neq  0$ and so $H_{s+1,2}\neq 0$.
\QED

The local cohomology modules $H_{s,t}$ are a good illustration for the failure of the methods of the previous section in the non-graded case.
For example, one cannot decide whether $H_{s,t}$ is zero just by looking at $F_{s,t}$: the vanishing of
$H_{s,t}$ depends on the characteristic of $C$!
Compare this situation to the following graded problem.

\begin{Theorem}[A Weaker Monomial Conjecture]
Let $T$ be a local ring with system of parameters $x_1, \dots, x_s$. For all $t\geq 0$ we have
$$\left(x_1 \cdot \ldots \cdot x_s\right)^t \notin \langle x_1^{s t}, \dots, x_s^{s t} \rangle.$$
\end{Theorem}

\demo
Let $S=\mathbb{Z}[A_1, \dots, A_s][X_1, \dots, X_s]$ where $\deg A_i=0$ and $\deg X_i=1$ for all $1\leq i\leq s$.
Following Hochster's argument we reduce to the problem of showing that
$$H^s_{\langle X_1, \dots, X_s \rangle}(S/fS)=0$$
where
$$f=(X_1 \cdot \ldots \cdot X_s)^t - \sum_{i=1}^s A_i X_i^{s t} .$$
Since $f$ is homogeneous and $\content(f)$ is the unit ideal, the result follows from Theorem \ref{theorem1}.
\QED


\begin{thebibliography}{KS}

\bibitem[BH]{BH} M.~Brodmann and M.~Hellus, \textit{Cohomological patterns of coherent sheaves over projective schemes,} J.~Pure Appl.~Algebra 172 (2002)\textbf{2-3}, pp. 165--182.

\bibitem[BS]{BS}M.~P.~Brodmann and R.~Y.~Sharp,
\textit{Local cohomology: an algebraic introduction with geometric applications,}
Cambridge University Press, 1998.

\bibitem[Ho]{Ho} M.~Hochster,
\emph{Canonical elements in local cohomology modules and the direct summand conjecture,}
J. Algebra \textbf{84} (2) (1983), pp.~503--553.

\bibitem[Hu]{Hu}C.~Huneke,
\textit{Problems on local cohomology},
in  \textit{Free resolutions in commutative algebra and algebraic geometry, Sundance 90},
ed.~D.~Eisenbud and C.~Huneke, Research Notes in Mathematics \textbf{2},
Jones and Bartlett Publishers, Boston, 1992, pp.~93--108.


\bibitem[K]{K} M.~Katzman,
\textit{An example of an infinite set of associated primes of a local cohomology module,}
J. Algebra \textbf{252}(1) (2002), pp.~161--166.


\bibitem[KS]{KS} M.~Katzman and R.~Y.~Sharp,
\emph{Some properties of top graded local cohomology modules,}
J. Algebra \textbf{259}(2) (2003) , pp.~599--612.

\bibitem[L]{L} G.~Lyubeznik,
private communication.

\bibitem[R]{R} P.~Roberts,
\emph{A computation of local cohomology,}
in
\emph{Commutative algebra: syzygies, multiplicities, and birational algebra (South Hadley, MA, 1992),}
Contemp. Math., \textbf{159},
Amer. Math. Soc., Providence, RI, 1994,
pp.~351--356.


\end{thebibliography}
\end{document}